\documentclass[final,5p,times,twocolumn]{elsarticle}

\usepackage{amssymb}
\usepackage{amsmath,amsfonts,mathrsfs}
\usepackage{lipsum}
\usepackage[T1]{fontenc}
\usepackage{newtxtext}
\usepackage{conjecture_perso}
\usepackage{algorithm}
\usepackage{algpseudocode}
\usepackage{graphicx}
\usepackage{textcomp}
\usepackage{xcolor}
\usepackage{hyperref}
\usepackage{csquotes}
\usepackage{tikz}
\usepackage{pgfplots}

\usetikzlibrary{arrows,shapes,positioning,shadows,trees,decorations.pathreplacing,calligraphy}
\usepgfplotslibrary{colormaps} 
\usepgfplotslibrary{groupplots}
\pgfplotsset{width=7cm,compat=1.18} 
\pgfplotsset{colormap={bwr}{[1cm]
        rgb255(0cm)=(0, 0, 255)
        rgb255(1cm)=(128, 128, 255)
        rgb255(2cm)=(255, 255, 255)
        rgb255(3cm)=(255, 128, 128)
        rgb255(4cm)=(255, 0, 0)}} 
        
\pgfplotsset{colormap={terrain}{[1cm]
        rgb255(0cm)=(51, 51, 153)
        rgb255(1cm)=(13, 125, 227)
        rgb255(2cm)=(0, 188, 148)
        rgb255(3cm)=(85, 221, 119)
        rgb255(4cm)=(197, 243, 141)
        rgb255(5cm)=(226, 217, 137)
        rgb255(6cm)=(170, 146, 107)
        rgb255(7cm)=(143, 111, 104)
        rgb255(8cm)=(199, 183, 179)
        rgb255(9cm)=(255, 255, 255)}} 

\newcommand{\bbR}{\mathbb{R}}

\newcommand{\bbE}{\mathbb{E}}

\newcommand{\tr}{\text{tr}}
\newcommand{\dd}{\mathrm{d}}
\newcommand{\W}{\mathrm{W}}
\newcommand{\tsp}{^\mathsf{T}}
\newcommand{\defeq}{\overset{\text{def}}{=}}

\newcommand{\muns}{\hat\mu^{\sigma,(n)}} 
\newcommand{\nuns}{\hat\nu^{\sigma,(n)}}

\newcommand{\Wns}{\mathbf{W}^{\sigma,(n)}}
\newcommand{\Wn}{\mathbf{W}^{(n)}}
\newcommand{\fns}{f^{\sigma,(n)}}
\newcommand{\fn}{f^{(n)}}
\newcommand{\Fnln}{F_{\mathcal{LN}}^{(n)}}
\newcommand{\Fn}{F^{(n)}}
\newcommand{\ftmin}{\textsc{ftmin} }
\newcommand{\gw}{\textsc{gw} }
\newcommand{\gpw}{\textsc{gpw} }
\newcommand{\wft}{\textsc{wft} }
\newcommand{\fisher}{\textsc{ftfish} }
\newcommand{\ks}{\textsc{ks} }
\newcommand{\kui}{\textsc{kui} }
\newcommand{\SW}{\mathrm{SW}}

\definecolor{primtext}{RGB}{43, 107, 204}

\journal{Journal of Statistical Planning and Inference}

\begin{document}

\begin{frontmatter}

\title{Two-sample test with Wasserstein distance on Gaussian samples based on a log-normal approximation}

\author[1]{Johann Clément-Cottuz}
\author[1]{Maxime Bérar} 
\author[2]{Gilles Gasso}

\affiliation[1]{organization={LITIS UR 4108, Univ Rouen Normandie}, 
            city={Rouen},
            postcode={76000}, 
            country={France}}
            
\affiliation[2]{organization={LITIS UR 4108, INSA Rouen Normandie}, 
            city={Rouen},
            postcode={76000}, 
            country={France}}

\begin{abstract}
In this article, we propose a two-sample test based on the Wasserstein $2$-distance between two 1D samples, where the first sample is assumed to be Gaussian based on an approximation of the underlying distribution by a log-normal distribution. A two-sample test is performed to decide whether two samples are drawn from the same distribution or not. We evaluate the proposed method on synthetic data and extend it to the multivariate case. We also show that this test is closely related to the Shapiro–Wilk normality test.
\end{abstract}


\end{frontmatter}

\section{Introduction}
Optimal transport is a well studied mathematical field (\cite{villani2003topics}, \cite{Peyre2020}) that defines the Wasserstein distance between two probability measures. This paper aims to use this metric to tackle the two-sample test problem, where one has to decide whether two samples are drawn from the same distribution. Wasserstein distance has been explored to address this problem (\cite{Matsui2016}, \cite{Ramdas2017}, \cite{Schefzik2021}), often with permutation tests (\cite{edgington2025}) to compute $p$-values (\cite{Matsui2016}, \cite{Hu2025}, \cite{Wang2022}, \cite{Schefzik2021}, \cite{Knijnenburg2009}). The Wasserstein distance has also been studied in the case of Gaussian samples by \cite{Rippl2016}, leading to the limit distribution of a Bures distance that approximates the Wasserstein distance between samples. However, the derived  two-sample test is not valid as the limit distribution advocated in \cite{Rippl2016} appears to be far from the true one in practice, even for large sample size. In this paper, we propose another characterization of the Wasserstein distance between Gaussian samples, which leads to a practical two-sample test. As in Mardia’s test (\cite{Mardia1970}), which relies on a fitted parametric distribution under the null hypothesis, we approximate the distribution of the Wasserstein distance between two Gaussian samples by fitting a log-normal distribution.

This paper is structured as follows. We recall some definitions and provide a brief overview of existing results on the Wasserstein distance between samples (Section \ref{sec:bib}). Then, we propose a fitted distribution for a Wasserstein-based statistic computed from two Gaussian samples, using Monte Carlo simulations to estimate suitable parameters (Section \ref{sec:conj}). Based on this distribution under the null hypothesis, we derive a two-sample test (Section \ref{sec:tstest}). We study the resulting test and compare it with other methods that do not rely on optimal transport, such as the classical tests reviewed in \cite{twosamplefrank}, which focus on Gaussian samples (Section \ref{sec:comp}). We also show that the proposed two-sample test can be viewed as a generalization of the Shapiro–Wilk normality test (Section \ref{sec:sw}). Finally, we extend the univariate test to the multivariate setting and illustrate its relevance to a change point detection problem (Section \ref{sec:gpw}).

\section{Two-sample test with a Gaussian sample}

\subsection{Wasserstein distance on one-dimensional distributions}
\label{sec:bib}

Let us consider the metric space $(\mathcal P(\bbR), \W_p)$ where $\W_p$ is the \textit{Wasserstein $p$-distance} defined for $p\ge 1$ by $$\W_p(\mu, \nu) = \left( \min_{\gamma \in \Gamma (\mu,\nu)} \int_{\bbR \times \bbR} | x-y|^p\, \dd \gamma(x,y)\right)^\frac 1p$$ where $\Gamma (\mu,\nu)$ is the set of probability measures over $\bbR\times \bbR$ with marginals $\mu$ and $\nu$.

Let define the empirical measure of a set of points $\{X_i\}_{1\le i\le n}\in \bbR$ by $\hat \mu^{(n)} \defeq \frac 1n \sum_{i=1}^n \delta_{X_i}$. In the literature, the semi-discrete optimal transport quantity $\W_2^2(\hat\mu^{(n)}, \mu)$, which measures the distance between the empirical measure $\hat\mu^{(n)}$ and its continuous counterpart $\mu$, has been extensively studied for a wide range of distributions with for instance, the derivation of sharp bounds 
for $\bbE (\W_2^2(\hat\mu^{(n)}, \mu))$ in 
 \cite{Bobkov2019} and \cite{Fournier2015}. In the Gaussian case, this expectation  has been investigated by \cite{Ledoux2019}, \cite{ledoux2019_2} and \cite{Berthet2020}. For example, \cite{Berthet2020} shows that $$\lim_{n\to +\infty} \frac{n}{\log\log n}\bbE \big(\W_2^2(\hat\mu^{(n)}, \mu)\big)=1.$$
 
Given another set of points $\{Y_i\}_{1\le i\le m}\in \bbR$ and its associated empirical measure $\hat\nu^{(m)}$, we can compare the two sets by evaluating the Wasserstein distance between their empirical measures, $\W_2(\hat\mu^{(n)},\hat\nu^{(m)})$. When both sets have the same size $n$, their squared distance can be expressed as a minimum over the set $S_n$ of all permutations of the indices in $\{1,\ldots, n\}$ (\cite{Peyre2020}): 
\begin{equation}
\label{eq:Wperm}
    \W_2^2(\hat \mu^{(n)},\hat\nu^{(n)}) = \min_{\tau\in S_n}\frac 1n\sum_{i=1}^n\left(X_{\tau (i)} - Y_i\right)^2.
\end{equation}

This formulation allows to compute distances between large sets of points efficiently, as the computation reduces to sorting the two samples. However, the sorting operation introduces the order statistics which considerably complicate the theoretical analysis of $\W_2^2(\hat \mu^{(n)},\hat\nu^{(n)})$. Nevertheless, there are few results on the two-sample case: as mentioned by \cite{delBarrio2019}, $\W_2^2(\hat\mu^{(n)}, \hat\nu^{ (n)})\to \W_2^2(\mu,\nu)$ almost surely. In particular, when $\mu = \nu$, we have that $\W_2^2(\hat\mu^{(n)}, \hat\nu^{(n)}) \to 0$ almost surely. The $\mu\ne \nu$ case is studied in \cite{Berthet2020_2}. Convergence results in the Wasserstein 1-distance case are given in \cite{Berthet2019}. In the Gaussian case, convergence results with two correlated samples are found in \cite{Berthet2020}. \cite{Rippl2016} proposes a method to circumvent the order statistics introduced by sorting: they estimate the Wasserstein 2-distance between two samples through the Bures distance (\cite{Bures1969}) between the two related normal distributions\footnote{Which is equal to their Wasserstein 2-distance (\cite{Givens1984})} $\mathcal N(m_X, \frac{n}{n-1}s_X^2)$ and $\mathcal N(m_Y,\frac{n}{n-1} s_Y^2)$, associated with the estimated means and variances of the samples, leading to $$\mathcal B(m_X, m_Y,\tfrac{n}{n-1}s_X^2, \tfrac{n}{n-1}s_Y^2) = (m_X - m_Y)^2 + \tfrac{n}{n-1}(s_X-s_Y)^2.$$ Here, $m_X = \frac 1n \sum_{i=1}^n X_i$ and $s_X^2 = \frac 1n\sum_{i=1}^n (X_i-m_X)^2$.

Based on that estimation, \cite{Rippl2016} derives the asymptotic distribution of $\widehat{\mathcal{GW}}_{n,n} \defeq \mathcal{B}(m_X, m_Y,\tfrac{n}{n-1}s_X^2, \tfrac{n}{n-1}s_Y^2)$ under the null hypothesis $\mathcal{N}(m_X,\tfrac{n}{n-1}s_X^2) = \mathcal{N}(m_Y,\tfrac{n}{n-1}s_Y^2)$. However, as shown by \cite{Irpino2015}, this quantity differs substantially from the Wasserstein 2-distance $\W_2^2(\hat \mu^{(n)},\hat\nu^{(n)})$. Indeed, define the correlation coefficient $$\rho_{X,Y} = \frac{\frac 1n \sum_{i=1}^n (X_{\tau (i)}-m_X)(Y_{i}-m_Y)}{\sqrt{\left(\frac 1n \sum_{i=1}^n (X_{\tau (i)}-m_X)^2 \right) \left(\frac 1n \sum_{i=1}^n (Y_{i}-m_Y)^2 \right) }},$$ we have 
\begin{align}
    \W_2^2(\hat\mu^{(n)}, \hat\nu^{(n)}) =\; & (m_X-m_Y)^2 + (s_X - s_Y)^2 \nonumber \\
    & \! +2s_Xs_Y(1-\rho_{X,Y}). \label{eq:irpino}
\end{align} In \cite{Irpino2015}, the terms are referred as \textit{mean}, \textit{size} and \textit{shape} respectively to characterize how the two distributions $\mu$ and $\nu$ may differ. We can express explicitly the relation between the empirical Wasserstein distance and $\widehat{\mathcal{GW}}_{n,n}$: $$\W_2^2(\hat\mu^{(n)}, \hat\nu^{(n)}) = \widehat{\mathcal{GW}}_{n,n} + 2s_Xs_Y(1-\rho_{X,Y}) - \tfrac{1}{n-1}(s_X - s_Y)^2.$$ The quantity $\widehat{\mathcal{GW}}_{n,n}$ does not include the $\textit{shape}$ term, that will allow us to distinguish $\mu$ from a non-normal distribution $\nu$ sharing its mean and variance. The shape term involves order statistics making the derivation of limit distribution of $\W_2^2(\hat\mu^{(n)}, \hat\nu^{(n)})$ a \textquote{difficult open problem} (\cite{Hallin2020}). Hence, in the sequel, we will propose a practical approximation of that limit distribution, and will apply it to the two-sample test problem.

\subsection{Proposition of a density function for $\W_2^2(\mu^{(n)}, \nu^{(n)})$ in the Gaussian case under null hypothesis}
\label{sec:conj}

Suppose that $\mu$ and $\nu$ are identical normal distributions $\mathcal N(\theta, \sigma^2)$ and that $\{X_i\}_{1\le i\le n}$ and $\{Y_i\}_{1\le i\le n}$ are two i.i.d. samples of size $n$ drawn from $\mu$. We denote by $\muns$ (resp. $\nuns$) the empirical measure associated with the points $\{X_i\}_{1\le i\le n}$ (resp. $\{Y_i\}_{1\le i\le n}$).

We are interested in the behavior of the random variable $\Wns \defeq \W_2^2(\muns, \nuns)$. The randomness of $\Wns$ arises from the randomness of the samples $X_i$ and $Y_i$. As a first property, note that the behavior of $\Wns$ does not depend on the mean of $\mu$. Therefore, we assume without loss of generality that the mean $\theta$ of $\mu$ is zero.

\begin{property}[Reduction to the $\theta=0$ case]
The distribution of $\Wns$ does not depend on the mean $\theta$ of the normal distribution $\mu$. 
\end{property}

\begin{proof}
Let $\{X_i\}_{1\le i\le n}$, $\{Y_i\}_{1\le i\le n}$ be independent random variables of law $\mathcal N (\theta,\sigma^2)$. We have

\begin{align*}
    \Wns &= \min_{\tau \in S_n} \frac 1n \sum_{i=1}^n \| X_{\tau (i)} - Y_i \|^2 \\
    &= \min_{\tau \in S_n} \frac 1n \sum_{i=1}^n \| (X_{\tau (i)}-\theta) - (Y_i-\theta) \|^2,
\end{align*}
the last expression being the Wasserstein distance between samples $\{X_i-\theta\}_{1\le i\le n}$, $\{Y_i-\theta\}_{1\le i\le n}$ of distribution $\mathcal N(0,\sigma^2)$.
\end{proof}

Another property that links $\Wns$ to $\mathbf{W}^{1,(n)}$, the variable in the case $\sigma = 1$, can be stated as:

\begin{property}[Reduction to the $\sigma = 1$ case]
\label{pp:reduc}
For any $n\ge 1$, $\Wns = \sigma^2\mathbf{W}^{1,(n)}$.
\end{property}

\begin{proof}
\begin{align*}
    \frac 1{\sigma^2} \W_2^2\left(\muns,\nuns\right) &= \frac{1}{\sigma^2}\min_{\tau\in S_n} \frac 1n\sum_{i=1}^n\left\|X_{\tau (i)} - Y_i\right\|^2 \\
    &= \min_{\tau\in S_n} \frac 1n\sum_{i=1}^n\left\|\frac{1}{\sigma}X_{\tau (i)} - \frac{1}{\sigma}Y_i\right\|^2 \\ &= \W_2^2\left(\hat\mu^{1,(n)},\hat\nu^{1,(n)}\right).
\end{align*}  
\end{proof}

Thus, the distribution of $\Wns$, denoted $\fns$, can be obtained from the distribution of $\mathbf{W}^{1,(n)}$ using the relation $\fns(t) = \frac{1}{\sigma^2} f^{1,(n)}\left(\frac{t}{\sigma^2}\right)$. Consequently, it suffices to study the standard normal case. We therefore write $\Wn$ for $\mathbf{W}^{1,(n)}$ and $\fn$ for $f^{1,(n)}$.

To properly characterize the behavior of $\Wn$ and to derive a statistical test based on this quantity, we require an explicit form of $\fn$ for any value of $n$. To our knowledge, no such results exist in the literature. Herein, we propose a conjecture for the approximation of $\fn$ for sufficiently large values of $n$, typically $n\ge 20$.

\begin{conjecture*}[Distribution of ${\Wn}$]
    For $n$ large enough, the distribution of $\Wn$ can be approximated by a log-normal distribution $\mathcal{LN}(\mu_{n},\tau_n)$.
\end{conjecture*}

\noindent
In order to apply this conjecture on $\Wns$, we use the following property:

\begin{property}[Generalization to any $\sigma$]
\label{pp:wns_distr}
If $\Wn\sim \mathcal{LN}(\mu_n,\tau_n)$, then $\Wns\sim \mathcal{LN}(\mu_{\sigma,n},\tau_n)$ where $\mu_{\sigma,n} = \mu_n+2\log(\sigma)$. 
\end{property}

\begin{proof}
Indeed, one can apply a logarithm on Property \ref{pp:reduc}, which directly yields the result.
\end{proof}

\noindent
Property \ref{pp:wns_distr} gives an explicit form for the density of $\Wns$: $$\fns(t) \approx \frac{1}{t\tau_n\sqrt{2\pi}}\exp\left( - \frac{[\log(t) - (2\log (\sigma) + \mu_{n})]^2}{2\tau_n^2} \right).$$ Hence, we can compute tables for $\tau_n$ and $\mu_{n}$ and deduce an approximate distribution of $\Wns$ for any $\sigma > 0$ and $n\ge 20$.

The cumulative distribution function (cdf) of $\Wns$ can then be approximated as $P\left(\Wns \leq t \right) \approx \Phi\left( \frac{\log(t) - \mu_{\sigma,n}}{\tau_n}\right)$ where $\Phi$ is the cdf of the standard normal distribution. Setting $t_{n,\alpha} = \exp(\tau_n z_{1-\alpha} + \mu_{\sigma,n})$ with $z_{1-\alpha} = \Phi^{-1}(1-\alpha)$ ensures that $P\left(\Wns\leq t_{n,\alpha}\right) \approx 1-\alpha$.\\

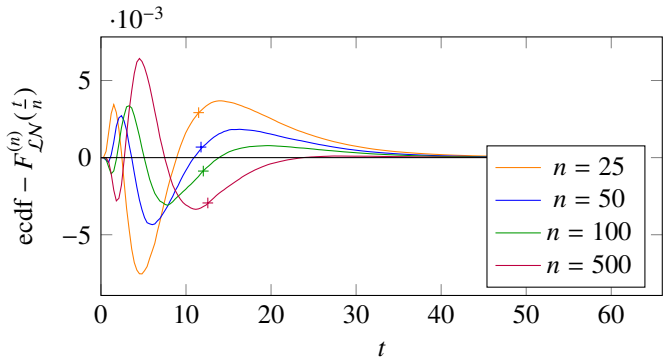
\begin{figure}[H]
\centering
\begin{tikzpicture}
    \begin{axis}[
      xlabel={$t$},
      ylabel={$\text{ecdf}-F_{\mathcal{LN}}^{(n)}(\frac{t}{n})$\vphantom{$\varepsilon_n([\Fnln]^{-1}(0.95))$}},
      grid=none,
      xmin=0,
      width=9cm,
      height=5cm,
      legend pos=south east
    ]
      \addplot[no markers, orange] table [col sep=space, header=true] {err_data25.txt};
      \addlegendentry{$n=25$}
      \addplot[no markers, blue] table [col sep=space, header=true] {err_data50.txt};
      \addlegendentry{$n=50$}
      \addplot[no markers, green!60!black] table [col sep=space, header=true] {err_data100.txt};
      \addlegendentry{$n=100$}
      \addplot[no markers, purple] table [col sep=space, header=true] {err_data500.txt};
      \addlegendentry{$n=500$}
      \addplot[black, domain=0:60] {0};
       \addplot [only marks, mark=+, color=orange] coordinates {(11.4863,0.0029165)};
      \addplot [only marks, mark=+, color=blue] coordinates {(11.76455, 0.0006835)};
        \addplot [only marks, mark=+, color=green!60!black] coordinates {(12.035,-0.000872)};
        \addplot [only marks, mark=+, color=purple] coordinates {(12.56552,-0.0029392)};
    \end{axis}
  \end{tikzpicture}
  \caption{Error between the estimated cumulative distribution of $n\Wn$ and estimated log-normal cdf for $n\in\{25, 50, 100, 500\}$.  The respective values at $t = [\Fnln]^{-1}(0.95)$ are displayed as $+$ marks}
\label{fig:graphe_ecdf_top}
\end{figure}

To gauge the numerical accuracy of this approximation, we compute the empirical cumulative distribution of $\Wn$ from $10^7$ realisations. We then compare it to the cumulative distribution function $\Fnln$ corresponding to the log-normal distribution $\mathcal{LN}(\mu_n,\tau_n)$, where $\mu_n$ and $\tau_n$ are estimated to fit the empirical distribution of $\Wn$. In figure \ref{fig:graphe_ecdf_top}, results are presented for the quantity $n\Wn$ \textendash\ that allows plotting on the same scale for different values of $n$ \textendash, and we observe that the difference is on the order of $10^{-3}$.

\begin{figure}[H]
\centering
\begin{tikzpicture}
      \begin{axis}[
        width=9cm,
        height=5cm,
        xlabel={$n$},
        ylabel={$\varepsilon_n([\Fnln]^{-1}(0.95))$},
        legend style={at={(1.05,1)}, anchor=north west},
        grid=none,
      ]
        \addplot[no markers, primbg!80!blue] table [x=x, y=y, col sep=space] {err95.txt};
        \addplot [only marks, mark=+, color=orange] coordinates {(25,0.0029165)};
        \addplot [only marks, mark=+, color=blue] coordinates {(50, 0.0006835)};
        \addplot [only marks, mark=+, color=green!60!black] coordinates {(100,-0.000872)};
        \addplot [only marks, mark=+, color=purple] coordinates {(500,-0.0029392)};
        \addplot[dotted, gray, domain=0:8000] {0};
      \end{axis}
    \end{tikzpicture}
  
\caption{Estimated value of $\varepsilon_n([\Fnln]^{-1}(0.95))$ in terms of $n$}
\label{fig:graphe_ecdf_bot}
\end{figure}
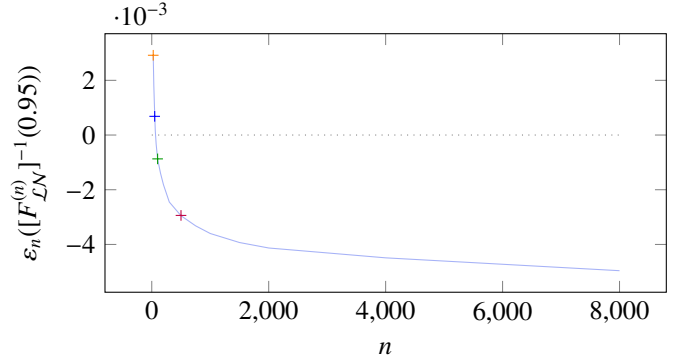

In order to assess the influence of the approximation on a test based on $\Fnln$, we define the error term $\varepsilon_n (t) = \Fn(t) - \Fnln(t)$, where $\Fn$ is the cumulative distribution function of $\Wn$. Figure \ref{fig:graphe_ecdf_bot} shows the values of $\varepsilon_n([\Fnln]^{-1}(0.95))$ in terms of $n$. For $20 \le n \le 500$, the error remains bounded by $3 \cdot 10^{-3}$. This finding is paramount in practice as $[\Fn]^{-1}(0.95)$ typically represents the threshold of the two-sample test under null hypothesis.

\subsection{Gaussian Wasserstein two-sample test}
\label{sec:tstest}
In the two-sample testing problem, we are given i.i.d.~samples $\{X_i\}_{1\le i\le n}$ and $\{Y_i\}_{1\le i\le m}$ drawn from unknown probability measures $\mu$ and $\nu$ respectively. A two-sample test is performed to decide whether to accept the null hypothesis $H_0$: $\mu = \nu$ or the alternative hypothesis $H_1$: $\mu \ne\nu$.\\

The results from Section \ref{sec:conj} allow us to establish the following test for one-dimensional samples. Assume that $\{X_i\}_{1\le i\le n} \sim \mathcal{N}(\theta, \sigma^2)$ are independent samples, and let $Y_1,\ldots, Y_n\in \bbR$. We propose the following test at significance level $\alpha$: 

\vspace{10pt}
\begin{center}
\noindent
\begin{minipage}{0.4\textwidth}
\begin{center}
    \textsc{Gaussian Wasserstein (gw) two-sample test}\\
\end{center}
\vspace{3pt}
\begin{itemize}[leftmargin=*]
\setlength\itemsep{0.5em}
\item[$\rhd$] $H_0$: points $\{Y_i\}_{1\le i\le n}$ are sampled independently from $\mathcal N(\theta,\sigma^2)$,
\item[$\rhd$] $H_1$: points $\{Y_i\}_{1\le i\le n}$ are not sampled independently from  $\mathcal N(\theta,\sigma^2)$.\\[-1pt]
\end{itemize}

\noindent
Reject $H_0$ if $$\W_2^2(\hat\mu_\sigma^{(n)}, \hat\nu^{(n)}) > \exp(\tau_n z_{1-\alpha} + \mu_{\sigma,n}).$$
\end{minipage}
\end{center}

\vspace{5pt}
\begin{center}
\rule{0.4\textwidth}{1pt}
\end{center}
\vspace{10pt}

Since the true value of $\sigma$ may be unknown, it can be estimated from the points $\{X_i\}_{1\le i\le n}$ using the method of \cite{Gurland1971} that amounts to modifying the denominator of the usual estimator of $\sigma$. Using this estimate $\hat \sigma$ in the test makes the significance level differ from the actual false alarm rate. In particular, applying the test at level $\alpha$ no longer guarantees a false alarm rate of $\alpha$. In a fixed $(\alpha, n)$ setting, we therefore perform a parameter study to determine an adjusted level $k^\gw_{\alpha,n}$. This quantity is defined such that, when used as the significance level in the \gw test with $\hat\sigma$, the resulting false alarm rate equals the target value $\alpha$. To do so, we reject $H_0$ if $\W_2^2(\hat\mu_\sigma^{(n)}, \hat\nu^{(n)}) > \exp(\tau_n z_{1-k^\gw_{\alpha,n}} + \mu_{\sigma,n})$.

\subsection{Comparison with standard two-sample tests}
\label{sec:comp}

The \gw test is designed to detect a shift from a normal distribution $\mathcal N (\theta,\sigma^2)$ to an alternative distribution. This change may take one of two forms: a shift in the parameters of a normal distribution, or a shift to a non-Gaussian distribution. We will examine both cases.

\subsubsection{The Gaussian case}

In a first place, let recall some standard tests for Gaussian samples. These include the $F$-test and the Student’s $t$-test. The $F$-test uses the unbiased variance estimates $\tfrac{n}{n-1}s_X^2$ and $\tfrac{m}{m-1}s_Y^2$, and evaluates the statistic $ Z = \frac{s_X^2}{s_Y^2}$ that follows a Fisher–Snedecor distribution of parameters $(n-1,m-1)$. The Student’s $t$-test computes the pooled variance estimator $S^2 = \frac{ns_X^2+ms_Y^2}{n+m-2}$ and hinges on the statistic $T = \sqrt{\frac{nm}{n+m}}\frac{m_X - m_Y}{S}$ that follows a Student's $t$ distribution with $n+m-2$ degrees of freedom.

A two-sample test for the equality of normal distributions is provided by \cite{ftmin}. Denoted as the \ftmin test, it combines the $p$-values of the $t$-test and the $F$-test by taking their minimum and rejects $H_0$ if this minimum is below $1-\sqrt{1-\alpha}$. Another standard method, coined the \fisher test, combines $p$-values according to \cite{fisher1932statistical} and rejects $H_0$ if $-2\log(TQ) \ge \chi^2_{4;1-\alpha}$, where $T$ and $Q$ are the $p$-values of the $t-$test and the $F$-test. A study by \cite{twosamplefrank} shows that the \ftmin test performs best when either mean or variance varies, whereas the \fisher test performs best when both mean and variance vary. Accordingly, we will compare the \gw test to these two tests.\\

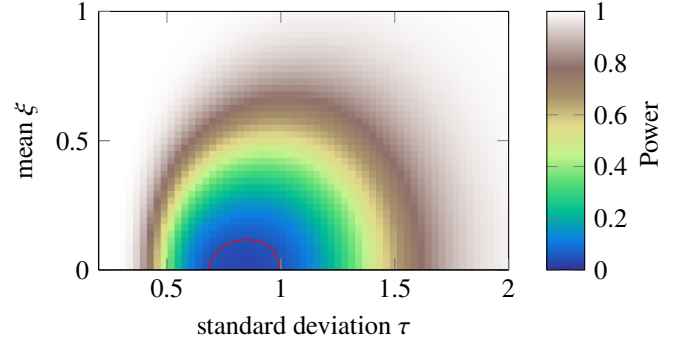
\begin{figure}[H]
\centering
\begin{tikzpicture}
\begin{groupplot}[
    group style={
        group size=2 by 1,
        horizontal sep=1.5cm
    },
    width=7cm,
    height=5cm,
    xlabel={standard deviation $\tau$},
    ylabel={mean $\xi$},
    colormap name=terrain,
    point meta min=0,
    point meta max=1
]

\nextgroupplot[
    title={}, view={0}{90}, colorbar, colorbar style={
    title=Power,
    title style={at={(2.7,0.3)}, anchor=west, rotate=90}
}
]
\addplot3 [
    mesh/ordering=y varies,
    surf,mesh/rows=40,
    shader=flat mean, colorbar
]
    table [x=x, y=y, z=z] {data_gw.txt};

\addplot3 [
    contour prepared,
    mesh/ordering=y varies,
    mesh/rows=40,
    mesh/cols=60,
    contour gnuplot={levels={0.05},draw color=red,labels=false},
]
    table [x=x, y=y, meta=z] {data_gw.txt};

\end{groupplot}
\end{tikzpicture}

\caption{Power of the \gw test with estimated $\sigma$ when the distribution of the first sample is $\mathcal N(0,1)$ and the second sample is $\mathcal N (\xi,\tau^2)$ with $n=40$ and $\alpha = 0.05$ (in red) and $k^\gw = \tfrac{5}{130}$}
\label{fig:test_power}
\end{figure}

We compute the statistical power of \gw on test cases of the form $X_i\sim \mathcal{N}(0,1)$ and $Y_i \sim\mathcal N (\xi,\tau^2)$ where the mean $\xi$ varies in $[0,1]$ and the standard deviation $\tau$ varies in $[0.2,2]$. More extreme values of $\xi$ and $\tau$ are not interesting as the test would almost always reject $H_0$. The parameter $\sigma$ is estimated from the samples $\{X_i\}_{1\le i\le n}$, and the resulting power of \gw is displayed as a contour plot in Figure \ref{fig:test_power}. One can see that \gw\!\!'s power gets below $\alpha$ when $\xi = 0$ and $\tau$ is slightly smaller than $1$ as denoted with a red line in the contour plot. 
 
 In this $n=40$ setup, the Wasserstein distance between a sample of variance $1$ and a sample of standard deviation $\tau$ chosen in this area tends to be smaller than the Wasserstein distance between two samples of variance $1$ at fixed mean. It can be intuited by observing that points from the sample of standard deviation $\tau$  have less extreme values than points from the sample of variance $1$.
 
\begin{figure}[H]
\centering
\begin{tikzpicture}
\begin{groupplot}[
    group style={
        group size=2 by 1,
        horizontal sep=0.2cm,
    },
    width=0.24\textwidth,
    height=4cm,
    xlabel={standard deviation $\tau$},
    ylabel={mean $\xi$},
]

\nextgroupplot[
    title={$P_\fisher - P_\gw$}, view={0}{90},
    colormap name=bwr,
    point meta min=-1,
    point meta max=1
]
\addplot3 [
    mesh/ordering=y varies,mesh/rows=10,
    surf, shader=interp ,
    colorbar
]
    table [x=x, y=y, z=z] {gw_fisher.txt};

\nextgroupplot[
    title={$P_\ftmin - P_\gw$},
    view={0}{90},  
    zlabel={}, colorbar, 
    colormap name=bwr,
    point meta min=-1,
    point meta max=1,
    ylabel={},
    xlabel={},
    ymajorticks=false,
    colorbar style={
    title=$\Delta$Power,
    title style={at={(2.3,0.2)}, anchor=west, rotate=90}
}
]
\addplot3 [
    mesh/ordering=y varies,mesh/rows=10,
    surf, shader=interp ,
    colorbar
]
    table [x=x, y=y, z=z] {gw_ftmin.txt};

\end{groupplot}
\end{tikzpicture}
\caption{Comparing \fisher\!\!'s and \ftmin\!\!'s power to \gw\!\!'s when the distribution of the first sample is $\mathcal N(0,1)$ and the second distribution is $\mathcal N (\xi,\tau^2)$ with $n=40$, $\alpha=0.05$ and $k^\gw = \tfrac{5}{130}$}
\label{fig:power_diff_gw_ft}
\end{figure}
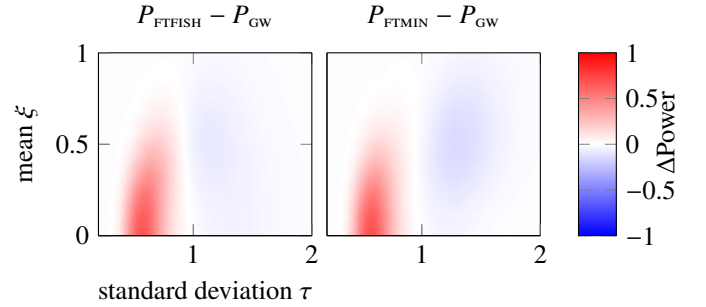

Within the same framework, we plot in Figure \ref{fig:power_diff_gw_ft} the difference in empirical power $P_\fisher - P_\gw$, between \fisher and \gw\!\!, with the red (resp. blue) regions indicating where \fisher (resp. \gw\!\!) performs best. We perform the same analysis for the \ftmin test. In both cases, a blue region appears when $\tau \ge \sigma$, indicating that the \gw test is particularly effective when handling increasing variances.

To address the inefficiency of the \gw test when the second sample has lower variance, we follow \cite{Heard2018} and construct a test that combines the \gw test with the $F$-test and the $t$-test. Specifically, we define the \wft test by combining the $p$-values of these three tests using Fisher’s method; that is, we reject $H_0$ if $$-2(\log p_{\gw} + \log p_{F} + \log p_t) \ge \chi_{6,1-\alpha}^2$$ where $p_\gw$, $p_F$ and $p_t$ are the $p$-values corresponding to the \gw test, $F$-test, and $t$-test, respectively. However, Fisher’s method requires the $p$-values to be independent, which is not the case here. Thus, the probability of falsely rejecting $H_0$ exceeds $\alpha$.

\begin{figure}[H]
    \centering
    \begin{tikzpicture}
  \begin{loglogaxis}[
    width=9cm,
    height=5cm,
    xlabel={$k^\wft$},
    ylabel={Empirical type-I error},
    title={},
    grid=major,
    legend style={at={(0.05,0.95)}, anchor=north west},
  ]
  
  \addplot [
      domain=0.00001:0.1,
      samples=100,
      thick,
      red,
    ]  {  0.5098 * x^(0.5582) };
    \addlegendentry{Interpolation: $\mathcal E_1 = a \cdot (k^\wft)^{b}$}

\addplot[
      only marks,
      mark=*,
      blue!80!black,
    ] table [col sep=space] {data_kwft.txt};

  \end{loglogaxis}
\end{tikzpicture}
    \caption{Type-I error of \wft for varying significance level $k^\wft$ at $n=40$}
    \label{fig:kwft}
\end{figure}
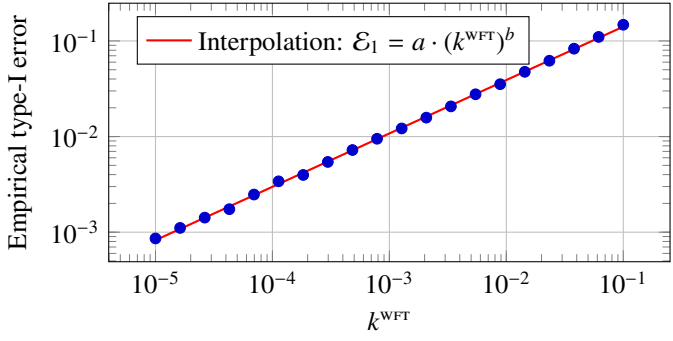

A suitable value for the provided significance level of the \wft test, denoted $k^\wft_{\alpha,n}$, can be obtained by observing that the type-I error $\mathcal E_1$ can be approximated by $\mathcal E_1 = a_n (k^\wft)^{b_n}$ where $a_n$ and $b_n$ depend only on $n$, as can be seen empirically in Figure \ref{fig:kwft}. Thus, $k^\wft_\alpha = (\alpha / a_n)^{1/b_n}$.

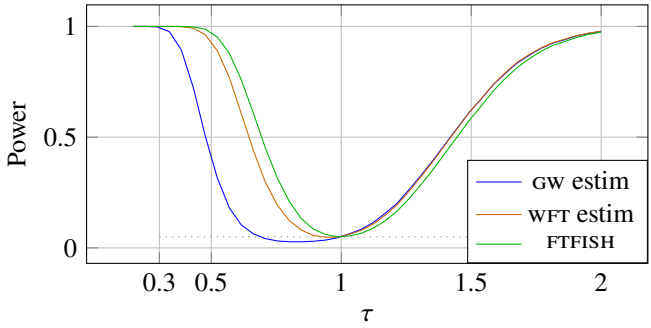
\begin{figure}[H]
    \centering
    \begin{tikzpicture}
      \begin{axis}[
        width=9cm,
        height=5cm,
        xlabel={$\tau$},
        xtick={0.3, 0.5, 1, 1.5, 2},
        ylabel={Power},
        legend style={at={(1,0)},anchor=south east},
        grid=major,
      ]
        \addplot[no markers, blue] table [x=tau, y=GW_estim, col sep=space] {power_wft_fish.txt};
        \addlegendentry{\gw estim}
    
        \addplot[no markers, orange!80!black] table [x=tau, y=GFT_estim, col sep=space] {power_wft_fish.txt};
        \addlegendentry{\wft estim}
    
        \addplot[no markers, green!70!black] table [x=tau, y=FTFISH, col sep=space] {power_wft_fish.txt};
        \addlegendentry{\fisher}
        
        \addplot[dotted, gray, domain=0.3:2] {0.05};
      \end{axis}
    \end{tikzpicture}
    \caption{Comparing \wft\!\!'s power to \gw\!\!'s and \fisher\!\!'s when the distribution of the first sample is $\mathcal N(0,1)$ and the second sample is $\mathcal N(0,\tau^2)$ with $n=40$, $\alpha = 0.05$, $k^\wft =\tfrac{5}{300}$ and $k^\gw  = \tfrac{5}{130}$}
    \label{fig:comp_var_gw_ftfish}
\end{figure}

In Figure \ref{fig:comp_var_gw_ftfish}, we plot the power of \fisher\!\!, \gw and \wft\!\!, with estimated $\sigma$ in terms of the second's sample standard deviation $\tau$.

We see in Figure \ref{fig:comp_var_gw_ftfish} that \wft is not as weak as \gw when $\tau \le \sigma$ and is still better than \fisher when $\tau \ge \sigma$. We can conclude that the $\wft$ test can be of interest in a fixed setup where $a_n$ and $b_n$ can be computed beforehand to obtain $k^\wft$. 

\subsubsection{The non-Gaussian case}

By construction, since the Wasserstein statistic accounts for both the mean and the variance, it is less sensitive to deviations captured by statistics that focus on only one of these features. However, it is not limited to comparisons between Gaussian samples: it can also distinguish between different distributions that share the same mean and variance, while \fisher and a test based on $\widehat{\mathcal{GW}}_{n,n}$ are unable to detect a difference in such a case. This robustness to non-normality in this setup is especially interesting in the case of change point detection where points from two distributions can be mixed (\cite{Truong2020}, Section 5.2.1). 

Table \ref{fig:nongaussian} reports the statistical power of \textsc{gw}, \textsc{wft}, \textsc{ftfish}, and other classical two-sample tests \textendash\ Kolmogorov–Smirnov (\textsc{ks}) (\cite{Pratt1981}) and Kuiper (\textsc{kui}) (\cite{Kuiper1960}) \textendash\ following test cases from (\cite{Rosenbaum2005}, Table 6) which also includes the cross-match (\textsc{cm}) test. We added other test cases, colored gray. We did not include a test based on $\widehat{\mathcal{GW}}_{n,n}$ because the limiting distribution provided in \cite{Rippl2016} deviates substantially from the empirical distribution for sample sizes $20 \le n \le 500$.

\begin{table}[H]
\centering
\begin{tabular}{lcccccc}
\hline
Distr. & \gw & \wft & \fisher & \ks & \kui & \textsc{cm}  \\
\hline
t$(1)$ & \textbf{1.00} & \textbf{1.00} & \textbf{1.00} & 0.12 & 0.57 & 0.19 \\
t$(2)$ & \textbf{0.95} & 0.94 & 0.92 & 0.05 & 0.19 & 0.08\\
$C_{0.25}^{0,100}$ & \textbf{1.00} & \textbf{1.00} & \textbf{1.00} & 0.08 & 0.38 & 0.19 \\
$C_{0.5}^{2,20}$ & \textbf{1.00} & \textbf{1.00} & \textbf{1.00} & 0.60 & 0.88 & 0.60\\
$C_{0.25}^{10,1000}$ & \textbf{1.00} & \textbf{1.00} & \textbf{1.00} & 0.14 & 0.48 & 0.19\\
$C_{0.5}^{10,1000}$ & \textbf{1.00} & \textbf{1.00} & \textbf{1.00} & 0.82 & \textbf{1.00} & 0.71\\
$\mathrm{E}(1) -1$ & 0.49 & 0.27 & 0.16 & 0.24 & \textbf{0.67} & 0.24 \\
\color{gray!80!black} $C_{0.25}^{1,5}$ &\color{gray!80!black} \textbf{0.74} &\color{gray!80!black} 0.71 & \color{gray!80!black}0.66 &\color{gray!80!black} 0.09 & \color{gray!80!black}0.14 &\color{gray!80!black} 0.06 \\
\color{gray!80!black} $\mathcal{L}(0,1)$ &\color{gray!80!black} \textbf{0.95} & \color{gray!80!black}\textbf{0.95} & \color{gray!80!black}0.94 & \color{gray!80!black}0.17 & \color{gray!80!black}0.60  &\color{gray!80!black} 0.14\\
\hline
\color{black} Type-I & 0.048 & 0.049  &  0.050 & 0.039 & 0.048 & 0.039 \\
\hline
\end{tabular}
\caption{Comparing two-samples tests on several distributions with $n=50$, $\alpha = 0.05$, $k^\gw = \frac{5}{130}$ and $k^\wft \approx 0.0159$. The first sample is drawn from $\mathcal N(0,1)$, the second is drawn from another distribution (first column), where t$(k)$ is the Student distribution with $k$ degrees of freedom, $C_{r}^{\xi, \tau^2}$ is a mixture of $\mathcal N(0,1)$ with probability $1-r$ and $\mathcal N(\xi,\tau^2)$ with probability $r$, $\mathrm{E}$ is the exponential distribution, and $\mathcal{L}$ is the logistic distribution}
\label{fig:nongaussian}
\end{table} 

In most cases tested here we observe that the \gw test obtains better rates.

\subsection{Test of normality}
\label{sec:sw}

A classical tool in data analysis that involves order statistics and Gaussian samples is the Shapiro–Wilk test (\cite{Shapiro1965}), which assesses the normality of a single sample. The test statistic is given by \begin{equation}
\label{eq:sw}
    \SW_n = \frac 1n \langle \mathbf{a}, \tilde{\mathbf{y}}\rangle^2
\end{equation} where $\tilde{\mathbf{y}}$ is the sorted vector of points $Y_i$, normalized by the sample standard deviation, and $\mathbf a$ is defined in \cite{Shapiro1965} as $\mathbf a = V^{-1} \mathbf{m} / \|V^{-1} \mathbf{m} \|$, where $\mathbf{m}$ is the vector of expected values of standard normal order statistics and $V$ the corresponding covariance matrix.

Following the spirit of the \gw test, we can construct a test of normality. We are given i.i.d. samples $\{Y_i\}_{1\le i\le n}$ drawn from an unknown probability distribution $\nu$. A test of normality is performed to decide whether to accept the null hypothesis $H_0$: $\nu$ is Gaussian or the alternative hypothesis $H_1$: $\nu$ is not Gaussian.

Instead of comparing two samples, we compare the sample vector $\mathbf{y}=(Y_1,\ldots, Y_n)$ to the reference vector $\mathbf a$.

Let us denote $\tilde{\mathbf y}=\frac{\mathbf y-\bar{\mathbf y}}{\sum_{i=1}^n (y_i-\bar{\mathbf y})^2}$ the reduced sample vector where $\bar{\mathbf y}$ is the sample mean. We reject $H_0$ when $\W^2_2(\mathbf a, \tilde{\mathbf y})$ is sufficiently large. In order to derive a decision threshold, we propose a conjecture analogous to that in Section \ref{sec:conj}.

\begin{conjecture*}[Distribution of ${\W^2_2(\mathbf a, \tilde{\mathbf y})}$ under $H_0$]

If the $Y_i$ are normally distributed and $n$ is sufficiently large, the distribution of $\W^2_2(\mathbf a, \tilde{\mathbf y})$ can be approximated by a log-normal distribution $\mathcal{LN}(\mu_{n},\tau_n)$.
\end{conjecture*}

\begin{figure}[H]
\centering
\begin{tikzpicture}
    \begin{axis}[
      xlabel={$x$},
      ylabel={Difference between cdfs},
      grid=none,
      xmin=0,
      width=9cm,
      height=5cm,
      legend pos=south east
    ]
      \addplot[no markers, orange] table [col sep=space, header=true] {err_data_a_30.txt};
      \addlegendentry{$n=30$}
\addplot[no markers, blue] table [col sep=space, header=true] {err_data_a_50.txt};
      \addlegendentry{$n=50$}
\addplot[no markers, green!60!black] table [col sep=space, header=true] {err_data_a_80.txt};
      \addlegendentry{$n=80$}
      \addplot[black, domain=0:0.2] {0};
    \end{axis}
  \end{tikzpicture}

\caption{Difference between the empirical cumulative distribution function and the associated log-normal cumulative distribution function}
\label{fig:ecdf_sw}
\end{figure}
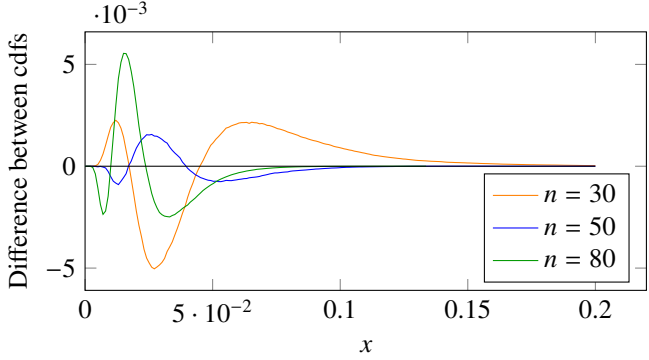

As in Figure \ref{fig:graphe_ecdf_top}, we assess the numerical accuracy of this approximation by comparing the empirical cumulative distribution of $n\W^2_2(\mathbf a, \tilde{\mathbf y})$ with  the corresponding log-normal cumulative distribution function. The results are gathered in Figure \ref{fig:ecdf_sw}. We observe that the approximation error is on the order of $10^{-3}$.

We perform the Wasserstein test of normality (\textsc{wn}) by rejecting $H_0$ if $$s_\alpha \le \W^2_2(\mathbf a, \tilde{\mathbf y})$$ where $s_\alpha = \exp (\tau_n \Phi^{-1}(1-\alpha)+\mu_n)$ and $\Phi$ is the cumulative distribution function of the standard normal. We can derive the Shapiro–Wilk test from the \textsc{wn} test:

\begin{property}[Deriving Shapiro–Wilk from \textsc{wn}]
\label{pp:sw}
For any $\alpha\in [0,1]$, $$s_\alpha \le \W^2_2(\mathbf a, \tilde{\mathbf y})\iff (1-\frac n2 s_\alpha)^2 \ge \SW_n$$ where $\SW_n$ has been defined is Equation \ref{eq:sw}.
\end{property}

\begin{proof}
Suppose the vector $\tilde{\mathbf y}$ sorted. Using that $\| \tilde{\mathbf y} \| =  \| \mathbf{a}\| = 1$, we have $$\W^2_2(\mathbf a, \tilde{\mathbf y})= \frac 1n \sum_{i=1}^n (\mathbf a_i - \tilde{\mathbf y}_i)^2 = \tfrac{2}{n}(1- \langle \mathbf a, \tilde{\mathbf y} \rangle) $$ where the first equality is given by (\cite{Peyre2020}, Remark 2.28).
Thus, since $\langle \mathbf a, \tilde{\mathbf y} \rangle$ is non-negative, $$s_\alpha \le \W^2_2(\mathbf a, \tilde{\mathbf y})\iff 1-\tfrac n2 s_\alpha\ge  \langle \mathbf a, \tilde{\mathbf y} \rangle \iff (1-\tfrac n2 s_\alpha)^2\ge  \langle \mathbf a, \tilde{\mathbf y} \rangle^2$$ which is the result.
\end{proof}

Property \ref{pp:sw} gives an approximation of the Shapiro–Wilk decision threshold for any value of $\alpha$, depending only on tables for $\mu_n$ and $\tau_n$. This threshold is usually approximated using \cite{Royston1995}, where the variable $\log (1-\SW_n)$ is normalized according to \cite{Royston1992}.  Table \ref{fig:sw_stat} compares the obtained threshold for $\alpha = 0.05$ and $\alpha = 0.02$ to the ones provided by \cite{Shapiro1965} and \cite{Hanusz2011}. We observe that using Property \ref{pp:sw} indeed provides similar thresholds as the ones given by the most recent method. This validates the use of our approximation for the Shapiro–Wilk test.

\begin{table}[H]
\centering

\begin{tabular}{c|cccc}
$n$ & 20 & 30 & 45 & 58\\ \hline
\textsc{wn} & 0.902 & 0.930 & 0.950 & 0.959\\
 \citeauthor{Hanusz2011}  & 0.904 & 0.931 & 0.949 & 0.959 \\
\citeauthor{Shapiro1965}  & 0.905 & 0.927 & 0.945 & $\times$ \\ \hline
\textsc{wn} & 0.880 & 0.914 &  0.939 & 0.951 \\
 \citeauthor{Hanusz2011}   & 0.884 & 0.916 & 0.938 & 0.951 \\
\citeauthor{Shapiro1965}  & 0.884 & 0.912 & 0.934 & $\times$ \\
\end{tabular}
  
\caption{Threshold of the Shapiro–Wilk at $\alpha = 0.05$ (top) and $\alpha = 0.02$ (bottom) derived from the \textsc{wn} one compared to the literature}
\label{fig:sw_stat}
\end{table}

\section{Deriving a $d$-dimensional  two-sample test}
\label{sec:gpw}

Herein, we lift the two-sample test in dimension $d\ge 2$. Suppose we have independent $d$-dimensional samples $X_1, \ldots, X_n \sim \mathcal{N}(\boldsymbol{\theta}, \Sigma)$.
Consider the orthonormal eigenvectors $\mathbf{v}_1, \ldots, \mathbf{v}_d$ of the covariance matrix $\Sigma$, with corresponding eigenvalues $\lambda_1, \ldots, \lambda_d$.
The projected samples $\mathbf{v}_j\tsp X_1,\ldots, \mathbf{v}_j\tsp X_n$ on vector $\mathbf{v}_j$ follow a normal distribution
$\mathcal{N}(\mathbf{v}_j\tsp \boldsymbol{\theta}, \mathbf{v}_j\tsp \Sigma \mathbf{v}_j)$. Thus, for each $1 \le j \le d$, we can perform a one-dimensional \gw test on the samples projected onto $\mathbf{v}_j$. We conduct the $d$ one-dimensional tests and reject $H_0$ if at least one of them rejects.

Regarding the second sample, let us assume that $Y_1, \ldots, Y_n$ follow a normal distribution $\mathcal N(\boldsymbol{\xi}, \Xi)$. By multiplying the points $X_i$ by $\tr(\Sigma)^{-\frac 12}$ and the points $Y_i$ by $\tr(\hat\Xi)^{-\frac 12}$, we restrict ourselves to the case $\tr(\Sigma ) = \tr(\Xi)$. The following Property~\ref{pp:dtest} ensures that under such hypothesis, the second sample has a larger variance in at least one of the 1D tests. Therefore, the \gw test is particularly well suited in this setting.

\begin{property}[Increasing variance in one projected sample]
\label{pp:dtest}
If $\tr(\Sigma) = \tr(\Xi)$, then there is an index $j$ such that $\mathbf{v}_j\tsp\Xi \mathbf{v}_j \ge \lambda_j$.
\end{property}

\begin{proof}
 Up to rotations, we can suppose that $(\mathbf{v}_j)_{1\le j\le d}$ is the canonical basis of $\bbR^d$. In that case, let $\mu_1,\ldots, \mu_d$ be the eigenvalues of $\Xi$ and $\mathbf{w}_1,\ldots, \mathbf{w}_d$ the corresponding eigenvectors. We can write $\mu_i = \sum_{k=1}^d \varepsilon_{k,i}\lambda_k$ with $\sum_{i=1}^d \varepsilon_{k,i} = 1$, the quantity $\varepsilon_{k,i}$ representing \textquote{the proportion of $\lambda_k$ going to $\mu_i$}. The variance of samples $\{Y_i\}_{1\le i\le n}$ projected on $\mathbf{v}_j = (0,\ldots, 1,\ldots , 0)$ is then $$\mathbf{v}_j\tsp \Xi \mathbf{v}_j = \sum_{i=1}^d \mu_i \mathbf{v}_j\tsp \mathbf{w}_i \mathbf{w}_i\tsp \mathbf{v}_j = \sum_{i=1}^d \mu_i \mathbf{w}_i[j]^2$$ where $\mathbf{w}_i[j]$ is the $j$-th coordinate of the vector $\mathbf{w}_i$. Now, the sum of the variances of projected samples over each vector $\mathbf{v}_1,\ldots, \mathbf{v}_d$ is 
\begin{align*}
\sum_{j=1}^d \mathbf{v}_j\tsp \Xi \mathbf{v}_j &= \sum_{j=1}^d\sum_{i=1}^d \left( \sum_{k=1}^d \varepsilon_{k,i}\lambda_k\right)\mathbf{w}_i[j]^2\\
&=\sum_{i=1}^d \left( \sum_{k=1}^d \varepsilon_{k,i}\lambda_k\right) \underbrace{\sum_{j=1}^d\mathbf{w}_i[j]^2}_{= 1}\\
&= \sum_{k=1}^d \lambda_k \underbrace{\sum_{i=1}^d \varepsilon_{k,i}}_{=1}=\sum_{k=1}^d \lambda_k.
\end{align*} Thus, there is an index $j$ such that $\mathbf{v}_j\tsp\Xi \mathbf{v}_j \ge \lambda_j$, otherwise the left hand side would be strictly smaller than the right hand side.
\end{proof}

In order to tackle the case where no increase of variance is observed despite a change in distribution, one can also compute the projected variance on additional random unit vectors. As mentioned in \cite{wckarl}, $n(n+1)/2$ vectors are sufficient to ensure a difference between the projected variances.

As we perform $d$ one-dimensional tests in order to achieve a $d$-dimensional test, we have to derive the level $\beta$ of the one-dimensional tests that leads to a $d$-dimensional test of level $\alpha$. When testing the projection on vector $\mathbf{v}_j$, we reject $H_0$ if an inequality of the form $\W_j > z_{j,\beta}$ is observed (see section \ref{sec:tstest}). As the projections $(\mathbf{v}_j\tsp X_i)_{1\le i\le n}, (\mathbf{v}_j\tsp Y_i)_{1\le i \le n}$ can be seen as coordinates in the basis $(\mathbf{v}_j)_{1\le j\le d}$, the events $\W_j > z_{j,\alpha}$ are independent with respect to $j$. Thus, the probability of rejecting $H_0$ in the $d$ dimensional test is \begin{equation}\label{eq:dtest_reject}
P\left(\bigvee_{j=1}^d (\W_j  > z_{j,\beta}) \right) = 1 - (1-\beta)^d.\end{equation} We then have to take $\beta = 1-(1-\alpha)^{\frac{1}{d}}$ so that $1 - (1-\beta)^d = \alpha$. Algorithm \ref{alg:gpw} summarizes the overall $d$-dimensional two-sample test.

\begin{algorithm}
    \caption{Gaussian Projected Wasserstein (\gpw\!\!) test}
    \label{alg:gpw}
    \begin{algorithmic}[1] 
        \Function{GPW}{$(X_i)_{1\le i\le n}$, $(Y_i)_{1\le i\le n}$, $\Sigma$, $\Xi$, $\alpha$} 
             \State $d\gets$ dimension of $X_1$
    		\State $\beta \gets 1 - (1-\alpha)^\frac 1d$ \Comment{Threshold of the 1-dimensional test}
    		\State $X_1,\ldots, X_n \gets \tr(\Sigma)^{-\frac 12}X_1,\ldots, \tr(\Sigma)^{-\frac 12}X_n$ \Comment{Trace normalization}
    		\State $Y_1,\ldots, Y_n \gets \tr(\Xi)^{-\frac 12}Y_1,\ldots, \tr(\Xi)^{-\frac 12}Y_n$
    \For{$\mathbf{v}\in \text{Sp}(\Sigma)$} \Comment{Iterates on eigenvectors of $\Sigma$}
    \State $\tilde X_1,\ldots,\tilde X_n \gets \mathbf{v}\tsp X_1,\ldots, \mathbf{v}\tsp X_n$ \Comment{Projection onto $v$}
    \State $\tilde Y_1,\ldots,\tilde Y_n \gets \mathbf{v}\tsp Y_1,\ldots, \mathbf{v}\tsp Y_n$
    \State $\hat\mu,\hat\nu \gets \frac 1n \sum_{i=1}^n \delta_{\tilde X_i}, \frac 1n \sum_{i=1}^n \delta_{\tilde Y_i}$
    \State $\sigma \gets \mathbf{v}\tsp \Sigma \mathbf{v}$
    \If{$\W^2_2 (\hat\mu,\hat\nu) > \exp (\tau_n z_\beta + \mu_{n,\sigma})$}
    \State \textbf{return} True
    \EndIf
    \EndFor
    \State \textbf{return} False
    \EndFunction
    \end{algorithmic}
\end{algorithm}

As the true values of $\Sigma$ and $\Xi$ are not known, they can be estimated with the points $X_i$ and $Y_i$ respectively. Note that under such an estimation of the covariance matrices, the error $\varepsilon_n([\Fnln ]^{-1}(\beta))$ is scaled by $d$ in the $d$-dimensional test as $1-(1-x)^d \sim_{x\to 0} dx$ in Equation \ref{eq:dtest_reject}.\\

We provide a change point detection (\cite{Truong2020}) application of the \gpw test, embedded in a window sliding method (\cite{Truong2020}, Section 5.2.1), in which the sample sequence is partitioned into small windows on which we use two-sample tests to assess the presence of a change point. When performing the \gpw test in this setup, we estimate $\Sigma$ and $\Xi$ online using a forgetting factor of value 0.998. We compare the \gpw test with a high-dimensional variant of the \fisher test, constructed the same way as in Algorithm \ref{alg:gpw}, except that there is no trace normalization. We also compare with \textsc{eagg} and \textsc{ediv} from \cite{Matteson2014} that use hierarchical clustering and permutation tests (\cite{edgington2025}). 

We use the following setup in $\bbR^8$: a sequence $(x_t)_{1\le t\le T}$ is generated with a true change point $0.3T \le t_c \le 0.9T$ such that points $(x_t)_{t \le t_c}$ are drawn from $\mu$ and points $(x_t)_{t_c < t}$ are drawn from $\nu$. The measures $\mu$ and $\nu$ are normally distributed, $\mu$ with mean $0$ and covariance diag(1,2,...,7,8), and $\nu$ with mean $(\delta m, \delta m, 0,\ldots, 0)$ and covariance diag($s(1),\ldots, s(8)$) where $s$ is a random permutation.

We give the average \textsc{RandIndex} (\cite{Truong2020}, Section 3.2.3) and the proportion of runs in which the \textsc{RandIndex} exceeds 0.9 for each method, across several values of $\delta m$, in Figure \ref{fig:cpd}. We also report the \textsc{RandIndex} of the \o\ method, which always returns that no change point has occurred. 

\begin{figure}[H]
\centering

\begin{tabular}{c|ccc|c}
$\delta m$ & 0 & 3 & 6 & Time\\ \hline
\gpw & \textbf{0.880 (0.62)} & \textbf{0.876} (0.61) & 0.884 (0.63)  & 0.11\\
\fisher & 0.858 (0.55) & 0.861 (0.56) & 0.869 (0.55) & 0.31\\
\textsc{eagg} & 0.425 (0.00) & 0.758 \textbf{(0.64)}& \textbf{0.998 (1.00)} & 2.09\\
\textsc{ediv} & 0.663 (0.10) & 0.664 (0.11) & 0.662 (0.10) & 4.96 \\
\o & 0.578 (0.00)& 0.581(0.00) & 0.580 (0.00) &  \\ 
\end{tabular}
  
\caption{Average \textsc{RandIndex}, proportion of runs in which the \textsc{RandIndex} exceeds 0.9 in parenthesis, and running time is seconds for several values of $\delta m$, with $T=2000$, $30$ windows, $\alpha = 0.025$}
\label{fig:cpd}
\end{figure}

We observe that \gpw outperforms \fisher in this setup, highlighting the strength of \gw when prior information about the first distribution is available. The change in mean affects only the performance of \textsc{eagg}, which achieves the best results when $\delta m = 6$, albeit at the cost of increased runtime.

\section{Conclusion}

We constructed a two-sample test based on the Wasserstein distance under the assumption that the first sample follows a Gaussian distribution. Our experiments show that this test is particularly powerful when the variance increases. Incorporating prior knowledge of the standard deviation of the first distribution can further improve the performance of the \gw test. This is especially relevant in change point detection settings, where the characteristics of a stationary regime are known.

We showed that using the Wasserstein distance for a normality test is equivalent to performing the Shapiro–Wilk test. Thus, we provided a new method to compute the Shapiro–Wilk decision threshold for any significance level $\alpha$. Among the natural extensions of this work, we could approximate the Wasserstein distance between samples drawn from other distributions, as has been done for the Shapiro–Wilk test (\cite{Shapiro72}). This approach would yield a family of two-sample tests based on a fitted log-normal distribution. A natural first choice would be the exponential distribution.

\section*{Acknowledgments}

This work was supported by the French Agence Nationale de
la Recherche (ANR), under grant ANR-23-CE23-0004 (project ODD).

\bibliographystyle{elsarticle-harv}
\biboptions{authoryear}

\bibliography{refs.bib}

\end{document}